\documentclass[11pt,reqno,oneside]{amsart}
\usepackage{amscd}
\usepackage{amsmath}
\usepackage{latexsym}
\usepackage{amsfonts}
\usepackage{amssymb}
\usepackage{amsthm}
\usepackage{graphicx}
\usepackage{verbatim}
\usepackage{mathrsfs}
\usepackage{enumerate}
\usepackage{hyperref}
\usepackage{fullpage}
\usepackage{xcolor}

\theoremstyle{plain}\newtheorem{definition}{Definition}[section]
\theoremstyle{plain}\newtheorem{theorem}{Theorem}[section]
\theoremstyle{plain}\newtheorem{lemma}[theorem]{Lemma}
\theoremstyle{plain}\newtheorem{coro}[theorem]{Corollary}
\theoremstyle{plain}
\theoremstyle{plain}

\numberwithin{equation}{section}

\newcommand{\norm}[1]{\left\|#1\right\|}

\newcommand{\R}{\mathbb{R}}
\newcommand{\be}{\begin{equation}}
\newcommand{\ee}{\end{equation}}
 \newcommand{\ba}{\begin{aligned}}
 \newcommand{\ea}{\end{aligned}}

  \newcommand{\ben}{\begin{enumerate}}
   \newcommand{\een}{\end{enumerate}}

\newcommand{\Rmnum}[1]{\expandafter\@slowromancap\romannumeral #1@}
\allowdisplaybreaks
\begin{document}

\title{ Local Well-posedness of Two Dimensional SQG Equation and Related Models}
\author[]{Huan Yu$^{1}$,\,Wanwan Zhang$^{2}$}

\address{$^{1}$School of Applied Science, Beijing Information Science and Technology University, Beijing  100192, P.R.China}
\email{ yuhuandreamer@163.com}

\address{$^{2}$School of Mathematical Sciences, Capital Normal University, Beijing, 100048, P.R.China}
\email{zhangwanwan153@163.com}


\begin{abstract}
In this paper, we present a new and elementary proof of  the local existence and uniqueness of
the classical solution to the Cauchy problem of the two-dimensional generalized surface quasi-geostrophic (SQG) equation  via the method of the contraction mapping principle. Also, same result holds true for a kind of transport equation with nonlocal velocity, of which finite time singularities were investigated in \cite{[D]} and \cite{[JZ]}.
\end{abstract}
\smallskip
\maketitle


\section{Introduction and main results}

In this paper, we consider   the Cauchy problem of  the  two-dimensional  transport equation with nonlocal velocity as follows
\begin{equation}\label{GSQG}
\left\{\ba
&\theta_{t}+u\cdot\nabla \theta = 0, ~(x,t)\in \R^{2}\times\R_{+},\\&u =\nabla^\bot(-\Delta)^{-1+\alpha}\theta,\\
&\theta(x,0)=\theta_{0}(x), \ea\ \right.
\end{equation}
where $0\leq\alpha\leq\frac{1}{2}$ and $\nabla^\bot=(\partial_2,-\partial_1).$
 The fractional Laplacian $(-\Delta)^s$ for any
$s\in\R$ is defined by the Fourier transform $\widehat{(-\Delta)^sf}(\xi)=|\xi|^{2s}\widehat{f}(\xi).$ Meanwhile, the expression of the second equation in \eqref{GSQG} implies that $u$ is divergence-free, that is, $\nabla\cdot u=\partial_{x_1}u_1+\partial_{x_2}u_2=0.$
Meanwhile, the second equation in \eqref{GSQG} can be written as
\begin{equation}\label{BS}
u(x,t)=C_\alpha P.V.\int_{\mathbb R^2}\frac{(x-y)^\bot}{|x-y|^{2+2\alpha}}\theta(y,t)dy,
\end{equation}
where $C_\alpha=-\frac{2\alpha\Gamma(\alpha)}{\pi2^{2-2\alpha}\Gamma(1-\alpha)}$.

When $\alpha=0$, \eqref{GSQG} is the two-dimensional incompressible  Euler equations, in which \eqref{BS} corresponds to the Biot-Savart law. The global regularity of the solution to the two-dimensional incompressible Euler equations  has been well-known (see  \cite{[Kato]}, \cite{[Koch]} and \cite{[CHEMIN]} and references therein).
When $\alpha=\frac12$, \eqref{GSQG}  reduces to the  $\text{SQG}$  equation which describes a famous approximation model of the nonhomogeneous fluid flow in a rapidly rotating 3D half-space (see \cite{[P]}).
The mathematical studies on the  SQG equation was first considered by Constantin, Majda and
Tabak in \cite{[CMT]} in which it was shown that the possible singularity of the SQG equation is an analogue to the 3D incompressible Euler equations. Also, as mentioned in \cite{[CMT]},  one can prove the local well-posedness of the  SQG equation in $H^k(\R^2)$ for some integers $k\geq3$ using the same techniques as for the incompressible Euler equation (see, e.g. \cite{[M]}). Recently, H. Inci \cite{[I]}  established the local well-posedness result for the  SQG equation in $H^s(\R^2)$ for $s>2$ by using a geometric approach. Existence of global weak solutions for the SQG equation was  proved by Resnick \cite{[R]}. The global regularity of smooth solution for the SQG equation has attracted considerable interests and progress has been made in this direction (see, e.g. \cite{[CCG]}, \cite{[Cor]}, \cite{[HeK]}). However, it is still  an outstanding open problem whether smooth solution  will blow up in finite time . 

A natural family of active scalar equations that interpolate between
the two-dimensional incompressible Euler equations and the SQG equation are given by \eqref{GSQG} with $0<\alpha<\frac{1}{2}.$
 This family has been called modified or generalized
SQG equation in the literature (see \cite{[KYZ]} and references therein).
 Kiselev, Yao, and Zlato\v{s} \cite{[KYZ]} established local $H^3$ patch solution  to the generalized SQG equation on the half-plane and  Kiselev, Ryzhik, Yao and Zlato\v{s}  \cite{[KRYZ]}  exhibited initial data
that lead to a singularity of the patch solution in finite time for $\alpha\in(0,\frac{1}{24})$.
Later, Gancedo and  Patel \cite{[GP]}  proved the local existence and uniqueness of  $H^2$ patch solution and also obtained the finite time singularity result in the half-space setting for $\alpha\in(0,\frac{1}{6})$.
In \cite{[YZJ]}, the authors proved that if the existence interval of the smooth solution to \eqref{GSQG} with $\alpha=\alpha_0\in[0,\frac{1}{2}]$ is $[0,T]$, then under the same initial data, the existence interval of the
solution to \eqref{GSQG} with $\alpha$ which is close to $\alpha_0$ will remain on $[0,T]$. As a direct byproduct, their results imply that the construction of the possible singularity of the smooth solution to
\eqref{GSQG} with $\alpha>0$ will be subtle, in comparison with the singularity presented in \cite{[KRYZ]}.

In this paper, we also consider the following transport equation with nonlocal velocity, in which the velocity is not  divergence-free:
\begin{equation}\label{tr-s-port}
\left\{\ba
&\theta_{t}+u\cdot\nabla \theta = 0, ~(x,t)\in \R^{2}\times\R_{+},\\&u =\nabla(-\Delta)^{-1+\alpha}\theta,\\
&\theta(x,0)=\theta_{0}(x), \ea\ \right.
\end{equation}
where $\nabla=(\partial_1, \partial_2).$
One  motivation of \eqref{tr-s-port} is that it can be regarded as a model equation for understanding the  $\text{SQG}$  equation.
Very recently, Dong and   Li \cite{[DL]} showed that  certain radial solutions
develop gradient blow-up in finite time to \eqref{tr-s-port} with $\alpha=\frac{1}{2}$. Formation of singularity for more general case was investigated in \cite{[D]} and \cite{[JZ]}. Chae in \cite{[C]} proved the local well-posedness of \eqref{tr-s-port} in N-dimensional space for $0<\alpha<1$. More precisely, the author obtained in \cite{[C]} the existence and uniqueness of local smooth solution to \eqref{tr-s-port} in $H^s(\R^N)$ with $s>\frac{N}{2}+1$ for $0<\alpha\leq\frac{1}{2}$ and  with $s>\frac{N}{2}+2$ for $\frac{1}{2}<\alpha\leq1$.

The aim of this paper is to
provide a new and elementary proof for  the local existence and uniqueness of
classical solutions of \eqref{GSQG} and \eqref{tr-s-port}  when $0<\alpha\leq\frac{1}{2}$ via  the contraction mapping principle. Our main results can be stated as

\begin{theorem}\label{the-1}Let $0<\alpha\leq\frac{1}{2}.$ For every $\theta_{0}\in H^{s}(\R^{2}),$ $s>1+2\alpha$, there exists a small time $T=T(\|\theta_{0}\|_{H^s(\R^{2})})$ such that $\eqref{GSQG}$ admits a unique solution $\theta\in C([0,T]; H^s(\R^{2})).$
\end{theorem}

\begin{theorem}\label{the-2}
Let $0<\alpha\leq\frac{1}{2}.$ For every $\theta_{0}\in H^{s}(\R^{2}),$ $s>1+2\alpha$, there exists a small time $T=T(\|\theta_{0}\|_{H^{s}(\R^{2})})$ such that $\eqref{tr-s-port}$ admits a unique solution $\theta\in C([0,T]; H^s(\R^2)).$
\end{theorem}
Let us describe the approach in some detail. Given $\theta_0\in H^s(\R^{2})$, we define a closed ball $B_T$ in the space $C([0,T];L^2(\R^2))$ (see \eqref{3.1} for the definition). Let $\theta\in B_T$, we take advantage of the global well-posedness result of the linear transport equation \eqref{3.2} to define a mapping $\mathcal{T}: B_T\rightarrow B_T$ (see \eqref{3.3}). Then we will prove that the mapping $\mathcal{T}$ has a unique fixed point $\theta$ in  $B_T$ through the contraction mapping principle, which is exactly the solution to $\eqref{GSQG}$ or $\eqref{tr-s-port}$. To this end, some delicate estimates on the singular integral involving the expression $u =\nabla^\bot(-\Delta)^{-1+\alpha}\theta$ in $\eqref{GSQG}$ or $u =\nabla(-\Delta)^{-1+\alpha}\theta$ in $\eqref{tr-s-port}$ will be applied. The uniqueness of the solutions to \eqref{GSQG} and \eqref{tr-s-port} is a consequence of the  $L^2$ energy estimate. In comparison with  Inci \cite{[I]} in which the case $\alpha=\frac12$ was treated, our proof is more elementary and our result holds true for $0<\alpha\le \frac 12$ in Theorem \ref{the-1}.
In comparison with Chae \cite{[C]} in which
the initial data belongs to $H^s(\R^2)$ with $s>2$ for $0<\alpha<\frac{1}{2}$, less regularity is imposed on the initial data in Theorem \ref{the-2}.  We can refer to \cite{[CCCG]} and \cite{[HKZ]} for the local well-posedness for the generalized SQG equation with more singular velocities, that is, with $\frac{1}{2}<\alpha<1$ in \eqref{GSQG}.
However, it is unclear that whether the current method of the contraction mapping principle can be applied to this more singular case.

The paper unfolds as follows.  In Section $2$, we recall some fundamental facts used in this paper. The proof of Theorem \ref{the-1} is given in Sections $3$. Sections $4$ is devoted to the proof of Theorem \ref{the-2}. In the end, some basic facts on Littlewood-Paley decomposition and Besov spaces are presented in Appendix A.

Throughout the paper, we denote $C$ an absolute constant which may depend on $\alpha>0$ and may be different from line to line. To emphasize the dependence on $\alpha$ we sometimes write $C_{\alpha}$.

\section{Preliminaries}
This section is concerned with some basic lemmas needed later. We begin with the well-known contraction mapping principle.
\subsection{The contraction mapping principle}
A mapping $\mathcal{T}$ from a normed linear space $\mathcal{V}$ into itself is called a contraction mapping if there exists a number $\kappa<1$ such that $\|\mathcal{T}x-\mathcal{T}y\|\leq\kappa\|x-y\|$ for all $x,y\in\mathcal{V}$.

\begin{lemma}\label{contration} {\rm (\cite{[GT]})} It holds that

{\rm(1)} A contraction mapping $\mathcal{T}$ in a Banach space $\mathcal{B}$ has a unique fixed point, that is there exists a unique solution $x\in\mathcal{B}$ of the equation $\mathcal{T}x=x$;

{\rm(2)} The result of (1) remains true if  the space $\mathcal{B}$ is replaced by any closed nonempty subset of $\mathcal{B}$.
\end{lemma}
\subsection{Hardy-Littlewood-Sobolev inequality}
In this subsection, we recall  the well-known Hardy-Littlewood-Sobolev inequality of fractional integral.

Let $0<\beta<n.$ The Riesz potential $I_{\beta}=(-\Delta)^{-\frac{\beta}{2}}$ is defined by
$$ I_\beta (f)(x)=\frac{1}{\gamma(\beta)}P.V.\int_{\mathbb{R}^n}\frac{f(y)}{|x-y|^{n-\beta}}dy$$ with $$\gamma(\beta)=\pi^{\frac{n}{2}}2^{\beta}\frac{\Gamma(\frac{\beta}{2})}{\Gamma(\frac{n}{2}-\frac{\beta}{2})},$$
where $\Gamma$ is the standard Gamma function.

\begin{lemma}{\rm(\cite{[S]})}\label{Hardy-Littlewood-Sobolev}
Let $0<\beta<n,$ $1<p<q<\infty,$ $\frac{1}{q}=\frac{1}{p}-\frac{\beta}{n}.$ Then, there exists a constant $C_{p,q}$ depending on $p$, $q$ such that $$\|I_\beta (f)\|_{L^{q}(\R^n)}\leq C_{p,q}\|f\|_{L^{p}(\R^n)}.$$
\end{lemma}
Taking $n=2$, $\beta=1-2\alpha$, $p=2$ and $q=\frac{1}{\alpha}$ in Lemma \ref{Hardy-Littlewood-Sobolev}, we obtain the following:
\begin{coro}\label{coro 2.3}
Let $0<\alpha<\frac{1}{2}$, for any $f\in L^2(\R^2)$, the following inequality holds true:
$$\|\Lambda^{-(1-2\alpha)}f\|_{L^{\frac{1}{\alpha}}(\R^2)}\leq C_\alpha\|f\|_{L^2(\R^2)},$$
where we have used the standard notation $\Lambda:=(-\Delta)^{\frac{1}{2}}.$
\end{coro}
\subsection{Existence of solutions for the transport equation}
In this subsection, we recall the well-posedness for the following linear transport equations:
\begin{equation}\label{2.1}
\left\{\ba
&\partial_tf+v\cdot\nabla f= g, ~(x,t)\in \R^{d}\times\R_{+},\\
&f_{|t=0}=f_{0}, \ea\ \right.
\end{equation}
where the functions $v:\R\times\R^d\rightarrow\R^d,~f_{0}:\R^d\rightarrow\R^N,$ and $g:\R\times\R^d\rightarrow\R^N$ are given.
For convenience, we present the definition and some properties of Besov spaces $B^\sigma_{p,r}(\R^d)$ in the Appendix A.
\begin{lemma}{\rm\cite{[BCD]}}\label{existence of transport}
Let $1\leq p \leq p_1\leq\infty$, $1<r<\infty$ and $p^\prime:=\frac{p}{p-1}$. Assume that
$$\sigma>-d\min(\frac{1}{p_1},\frac{1}{p^\prime})~or~\sigma>-1-d\min(\frac{1}{p_1},\frac{1}{p^\prime})~~~if~~{\rm div} v=0.$$
Let $f_{0}\in B^\sigma_{p,r}(\R^d)$, $g\in L^1([0,T];B^\sigma_{p,r}(\R^d))$ and $v$ be a time-dependent vector field such that $v\in L^\rho([0,T];B^{-\lambda}_{\infty,\infty}(\R^d))$ for some $\rho>1$ and $\lambda>0,$ and $\nabla v\in L^1([0,T];B^{\sigma-1}_{p_1,r}(\R^d))$ with $\sigma>1+\frac{d}{p_1}$.

Then the equation \eqref{2.1} has a unique solution $f$ in $C([0,T];B^\sigma_{p,r}(\R^d))$.
Moreover, the following estimate holds true: for a.e. $t\in[0,T],$
\begin{eqnarray}\begin{split}\label{2.2}
&\|f\|_{\widetilde{L}^\infty_t(B^\sigma_{p,r}(\R^d))}\\
\leq&\Big(\|f_0\|_{B^\sigma_{p,r}(\R^d)}+\int_0^te^{-C\int_0^{t^\prime}V_{p_1}^\prime(\tau)d\tau}\|g(t^\prime)\|_{B^\sigma_{p,r}(\R^d)}dt^\prime\Big)
e^{C\int_0^tV_{p_1}^\prime(t^\prime)dt^\prime},
\end{split}
\end{eqnarray}
with $V_{p_1}^\prime(t)=:\|\nabla v(t)\|_{B^{\sigma-1}_{p_1,r}(\R^d)}$.
\end{lemma}
Consequently, taking $d=2$, $N=1$, $f=\omega$, $v=u$ and $g=0$ in Lemma \ref{existence of transport} immediately implies the existence and uniquess of the solution to the following scalar transport equation
\begin{equation}\label{2.3}
\left\{\ba
&\partial_t\omega+u\cdot\nabla \omega=0,\\
&\omega(\cdot,0)=\theta_{0}, \ea\ \right.
\end{equation}

\begin{coro}\label{existence of scalar transport}
Let $\alpha\in(0,\frac{1}{2}]$ and $\theta_0\in H^s(\R^2),$ $s>1+2\alpha$. For $\theta\in L^\infty(0,T;H^s(\R^2))$, define $u=\nabla^\bot\Lambda^{-2+2\alpha}\theta.$
Then there exists a unique solution $\omega\in C([0,T]; H^s(\R^2))$ to \eqref{2.3}, satisfying
\begin{eqnarray}\label{2.4}
\|\omega\|_{L^\infty(0,T;H^s(\R^2))}\leq\|\theta_0\|_{H^s(\R^2)}e^{CT\|\theta\|_{L^\infty(0,T; H^s(\R^2))}}.
\end{eqnarray}
\end{coro}
{\bf Proof.} \textbf{Case 1: $\alpha=\frac{1}{2}$.} In this case, $u=\nabla^\bot\Lambda^{-1}\theta=\mathcal{R}^\perp\theta$, where $\mathcal{R}$ denotes the Riesz transform.
We take $p=p_1=2$, $r=2$ and $\sigma=s>2$ in Lemma \ref{existence of transport}.
Since $\theta\in L^\infty(0,T;H^s(\R^2))$, we have $u=\mathcal{R}^\perp\theta$ is also in $L^\infty(0,T;H^s(\R^2)).$ Therefore we have $\nabla u\in L^1([0,T];B^{s-1}_{2,2}(\R^2))$. Furthermore, by the standard embeddings of Besov spaces (see Lemma \ref{embeddings} in the Appendix A), we have $$H^s(\R^2)=B^s_{2,2}(\R^2)\hookrightarrow B^{s-1}_{\infty,2}(\R^2)\hookrightarrow B^{s-1}_{\infty,\infty}(\R^2)\hookrightarrow B^{-\lambda}_{\infty,\infty}(\R^2),$$
for any $\lambda>0$, since $s>2$.
Therefore, we obtain $u\in L^\rho([0,T];B^{-\lambda}_{\infty,\infty}(\R^2))$ for any $\rho>1$ and $\lambda>0.$
Then, by virtue of Lemma \ref{existence of transport}, we obtain that there exists a unique solution $\omega\in C([0,T]; H^s(\R^2))$ to \eqref{2.3}. Furthermore, it follows from \eqref{2.2} that, for a.e. $t\in[0,T]$,
\begin{eqnarray*}
\|\omega(t)\|_{H^s(\R^2))}&=&\|\omega(t)\|_{B^s_{2,2}(\R^2)}\\
&\leq&\|\theta_0\|_{B^s_{2,2}(\R^2)}e^{C\int_0^t\|\nabla u(\tau)\|_{B^{s-1}_{2,2}(\R^2)}d\tau}\\
&\leq&\|\theta_0\|_{H^s(\R^2)}e^{C\int_0^t\|\nabla u(\tau)\|_{H^{s-1}(\R^2)}d\tau}\\
&\leq&\|\theta_0\|_{H^s(\R^2)}e^{C\int_0^t\|u(\tau)\|_{H^s(\R^2)}d\tau}\\
&\leq&\|\theta_0\|_{H^s(\R^2)}e^{C\int_0^t\|\theta(\tau)\|_{H^s(\R^2)}d\tau}\\
&\leq&\|\theta_0\|_{H^s(\R^2)}e^{CT\|\theta\|_{L^\infty(0,T; H^s(\R^2)}},
\end{eqnarray*}
where we have used the equivalence between the Besov space $B^s_{2,2}(\R^2)$ and the standard sobolev space  $H^s(\R^2)$ and the fact that the Riesz transform $\mathcal{R}$ is bounded on the Sobolev space $H^s(\R^2).$
Therefore, we obtain \eqref{2.4}.

\textbf{Case 2: $0<\alpha<\frac{1}{2}$.} This case is more involved. $u=\nabla^\bot\Lambda^{-2+2\alpha}\theta=\Lambda^{-(1-2\alpha)}\mathcal{R}^\perp\theta.$
Similarly, we take $p=2$, $p_1=\frac{1}{\alpha}$, $r=2$ and $\sigma=s>1+2\alpha$ in Lemma \ref{existence of transport}. Since $\theta\in L^\infty(0,T;H^s(\R^2))$, we have $u\in L^\infty(0,T;B^{s}_{\frac{1}{\alpha},2}(\R^2))$. Indeed, by the definition of the norm of Besov space $B^{s}_{\frac{1}{\alpha},2}(\R^2)$ and Corollary  \ref{coro 2.3}, we have
\begin{eqnarray}\label{2.5}
\|u\|_{B^{s}_{\frac{1}{\alpha},2}(\R^2)}&=&\Big(\displaystyle{\sum_{j\geq-1}}2^{2sj}\norm{\Delta_{j} u}_{L^{\frac{1}{\alpha}}(\R^2)}^{2}\Big)^{\frac{1}{2}}\nonumber\\
&=&\Big(\displaystyle{\sum_{j\geq-1}}2^{2sj}\norm{\Delta_{j} \Lambda^{-(1-2\alpha)}\mathcal{R}^\perp\theta}_{L^{\frac{1}{\alpha}}(\R^2)}^{2}\Big)^{\frac{1}{2}}\nonumber\\
&=&\Big(\displaystyle{\sum_{j\geq-1}}2^{2sj}\norm{\Lambda^{-(1-2\alpha)}\Delta_{j} \mathcal{R}^\perp\theta}_{L^{\frac{1}{\alpha}}(\R^2)}^{2}\Big)^{\frac{1}{2}}\nonumber\\
&\leq& C\Big(\displaystyle{\sum_{j\geq-1}}2^{2sj}\norm{\Delta_{j} \mathcal{R}^\perp\theta}_{L^2(\R^2)}^{2}\Big)^{\frac{1}{2}}\nonumber\\
&\leq& C\Big(\displaystyle{\sum_{j\geq-1}}2^{sj}\norm{\Delta_{j}\theta}_{L^2(\R^2)}^{2}\Big)^{\frac{1}{2}}\nonumber\\
&=&C\|\theta\|_{B^s_{2,2}(\R^2)}\nonumber\\
&=&C\|\theta\|_{H^s(\R^2)}.
\end{eqnarray}
Therefore, we obtain $u\in L^\infty(0,T;B^{s}_{\frac{1}{\alpha},2}(\R^2))$.
Similarly, by the standard embeddings of Besov spaces, we have $$B^s_{\frac{1}{\alpha},2}(\R^2)\hookrightarrow B^{s-2\alpha}_{\infty,2}(\R^2)\hookrightarrow B^{s-2\alpha}_{\infty,\infty}(\R^2)\hookrightarrow B^{-\lambda}_{\infty,\infty}(\R^2),$$
for any $\lambda>0$, since $s>1+2\alpha$.
Therefore, we obtain $u\in L^\rho([0,T];B^{-\lambda}_{\infty,\infty}(\R^2))$ for any $\rho>1$ and $\lambda>0.$
By virtue of $u\in L^\infty(0,T;B^{s}_{\frac{1}{\alpha},2}(\R^2))$ and Bernstein inequality (see Lemma \ref{bern} in the Appendix A), we can obtain $\nabla u\in L^1([0,T];B^{s-1}_{\frac{1}{\alpha},2}(\R^2))$.
Then, by virtue of Lemma \ref{existence of transport}, we prove that there exists a unique solution $\omega\in C([0,T]; H^s(\R^2))$ to \eqref{2.3}. Furthermore, it follows from \eqref{2.2}, Bernstein inequality and \eqref{2.5} that, for a.e. $t\in[0,T]$,
\begin{eqnarray*}
\|\omega(t)\|_{H^s(\R^2))}&=&\|\omega(t)\|_{B^s_{2,2}(\R^2)}\\
&\leq&\|\theta_0\|_{H^s(\R^2)}e^{C\int^t_0\|\nabla u(\tau)\|_{B^{s-1}_{\frac{1}{\alpha},2}(\R^2)}d\tau}\\
&\leq&\|\theta_0\|_{H^s(\R^2)}e^{C\int^t_0\|u(\tau)\|_{B^{s}_{\frac{1}{\alpha},2}(\R^2)}d\tau}\\
&\leq&\|\theta_0\|_{H^s(\R^2)}e^{C\int^t_0\|\theta(\tau)\|_{H^s(\R^2)}d\tau}\\
&\leq&\|\theta_0\|_{H^s(\R^2)}e^{CT\|\theta\|_{L^\infty(0,T; H^s(\R^2)}},
\end{eqnarray*}
which implies \eqref{2.4}.
The proof of Corollary \ref{existence of scalar transport} is complete.

\section{Proof of Theorem \ref{the-1} }
In this section, we will prove Theorem \ref{the-1} by using the contraction mapping principle.
Before that, we introduce some notations. Given $s>1+2\alpha$ , set
\begin{equation}\begin{split}\label{3.1}
B_T=&\{\theta\in L^\infty(0,T;H^s(\R^2))\cap C([0,T];L^2(\R^2)):
ess\sup_{0\leq t\leq T} \|\theta(t)\|_{H^{s}(\R^2)}\leq M,\\&
 ~\theta(\cdot,0)=\theta_{0}\in H^s(\R^2)\},
 \end{split}
\end{equation}
where $M=2\|\theta_{0}\|_{H^{s}(\R^{2})}$ and $T>0$ is later to be determined.
We will construct the solution to \eqref{GSQG} as a fixed point of a mapping $\mathcal{T}:B_T\rightarrow B_T$.
For any $\theta\in B_T$,
let us consider the following transport equation
\begin{equation}\label{3.2}
\left\{\ba
&\partial_t\omega+u\cdot\nabla \omega=0,\\
&\omega(\cdot,0)=\theta_{0}, \ea\ \right.
\end{equation}
where $$u=\nabla^\bot\Lambda^{-2+2\alpha}\theta.$$
According to Corollary \ref{existence of scalar transport}, there exists a unique solution $\omega \in C([0,T];H^s(\R^2))$ to \eqref{3.2}. This process allows us to define the mapping
\begin{eqnarray}\label{3.3}
\mathcal{T}(\theta)(x,t)=\omega(x,t).
\end{eqnarray}
We claim that the mapping $\mathcal{T}:B_T\rightarrow B_T$ has exactly one fixed point in $B_T$. For clarity, we state it as the following key lemma, which will be used to prove  Theorem \ref{the-1}.
\begin{lemma}\label{existence of fixed-point}
The mapping $\mathcal{T}:B_T\rightarrow B_T$ defined by \eqref{3.3} has exactly one fixed point in $B_T$.
\end{lemma}
{\bf Proof.} The proof is divided into three steps.

{\it Step 1. $B_T$ is a closed nonempty subset of $C([0,T];L^2(\R^2))$.}

It is not difficult to prove that $B_T$ is a closed nonempty subset of $C([0,T];L^2(\R^2))$.
$B_T$ is nonempty, since $\theta_0\in B_T$.
To show that $B_T$ is closed in $C([0,T];L^2(\R^2))$, we assume that $\theta_n \in B_T$ and
\begin{eqnarray}\label{3.4}
\|\theta_n-\theta\|_{C([0,T];L^2(\R^2))}\rightarrow0,
\end{eqnarray}
as $n\rightarrow\infty$. We have to show that  $\theta\in B_T$.
Indeed,
 for a.e. $t\in [0,T],$ $\forall N>0$, by virtue of the triangle inequality and the Plancherel identity for the Fourier transform, we can arrive at
\begin{eqnarray}\label{3.5}
&&\Big(\int_{\mathbb B(0,N)}(1+|\xi|^2)^s|\widehat{\theta}(t,\xi)|^2 d\xi\Big)^\frac{1}{2} \nonumber\\
&\leq& \Big(\int_{\mathbb B(0,N)}(1+|\xi|^2)^s|\widehat{\theta}(t,\xi)-\widehat{\theta_n}(t,\xi)|^2 d\xi\Big)^\frac{1}{2}+\Big(\int_{\mathbb B(0,N)}(1+|\xi|^2)^s|\widehat{\theta_n}(t,\xi)|^2 d\xi\Big)^\frac{1}{2} \nonumber\\
&\leq& (1+N^2)^\frac{s}{2}\|\theta_n-\theta\|_{C([0,T];L^2(\R^2))}+\|\theta_n\|_{L^\infty(0,T;H^s(\R^2))} \nonumber\\
&\leq& (1+N^2)^\frac{s}{2}\|\theta_n-\theta\|_{C([0,T];L^2(\R^2))}+M.
\end{eqnarray}
Fix $N >0$, letting $n\rightarrow\infty$ in \eqref{3.5} and using \eqref{3.4}, we have
$$\Big(\int_{\mathbb B(0,N)}(1+|\xi|^2)^s|\widehat{\theta}(t,\xi)|^2 d\xi\Big)^\frac{1}{2}\leq M.$$
Since $N >0$ is arbitrary, we obtain $\theta\in L^\infty([0,T];H^s(\R^2))$ and $$\|\theta\|_{L^\infty(0,T;H^s(\R^2))}\leq M.$$
Moreover, it is clear that $\theta\in C([0,T];L^2(\R^2))$ and $\theta(\cdot,0)=\theta_{0}\in H^s(\R^2).$
Therefore,
\begin{eqnarray*}
\theta\in B_T,
\end{eqnarray*}
which shows that $B_T$ is a closed subset of $C([0,T];L^2(\R^2))$.\\[2mm]
{\it Step 2. $\mathcal{T}$ maps $B_T$ into $B_T$.}
In view of Corollary \ref{existence of scalar transport}, it is clear that $\omega\in C([0,T];L^2(\R^2))$ and $\omega(\cdot,0)=\theta_0\in H^s(\R^2).$
Besides, by \eqref{2.4}, we can obtain
\begin{eqnarray*}
\|\omega\|_{L^\infty(0,T;H^s(\R^2))}&\leq&\|\theta_0\|_{H^s(\R^2)}e^{CT\|\theta\|_{L^\infty(0,T; H^s(\R^2))}}\\
&\leq&\frac{1}{2}Me^{CMT}\\
&\leq&M,
\end{eqnarray*}
provided $T$ is chosen sufficiently small.
This shows that $\omega\in B_T.$
Therefore, $\mathcal{T}$ maps $B_T$ into $B_T$.\\[2mm]
{\it Step 3. The mapping $\mathcal{T}:B_T\rightarrow B_T$ is strictly contractive in the topology of $C([0,T];L^2(\R^2))$.}

 Suppose that $\theta_i\in B_T, u_i=\nabla^\perp\Lambda^{-2+2\alpha}\theta_i$ and $\omega_i=\mathcal{T}\theta_i(i=1,2).$
It follows from \eqref{3.2} that
\begin{equation}\label{3.6}
\left\{\ba
&\partial_t\omega_1+u_1\cdot\nabla\omega_1=0,\\
&\partial_t\omega_2+u_2\cdot\nabla\omega_2=0.\ea\ \right.
\end{equation}
Subtracting the equations in \eqref{3.6}, we get
\begin{eqnarray}\label{3.7}
\partial_t(\omega_1-\omega_2)+(u_1-u_2)\cdot\nabla\omega_1+u_2\cdot\nabla(\omega_1-\omega_2)=0.
\end{eqnarray}
Taking the inner product of \eqref{3.7} with $\omega_1-\omega_2$ in $L^2(\R^2)$ and using the divergence-free condition,
integration by parts, H\"{o}lder inequality, Lemma \ref{CZ inequality} (used when $\alpha=\frac{1}{2}$), Corollary \ref{coro 2.3} (used when $\alpha\in(0,\frac{1}{2})$) and the standard Sobolev embedding, we  have
\begin{eqnarray*}
\begin{split}
&\frac{1}{2}\frac{d}{dt}\|\omega_1(t)-\omega_2(t)\|^2_{L^2(\R^2)}\\=&-\int_{\mathbb{R}^2}((u_1-u_2)\cdot\nabla\omega_1)(\omega_1-\omega_2) dx-\int_{\mathbb{R}^2}(u_2\cdot\nabla(\omega_1-\omega_2))(\omega_1-\omega_2)dx\\
=&-\int_{\mathbb{R}^2}((u_1-u_2)\cdot\nabla\omega_1)(\omega_1-\omega_2)dx\\
\leq&\|u_1-u_2\|_{L^{\frac1\alpha}(\R^2)}\|\nabla\omega_1\|_{L^{\frac{2}{1-2\alpha}}(\R^2)}\|\omega_1-\omega_2\|_{L^2(\R^2)}\\
\leq&C\|\theta_1-\theta_2\|_{L^2(\R^2)}\|\omega_1\|_{H^s(\R^2)}\|\omega_1-\omega_2\|_{L^2(\R^2)}\\
\leq& CM\|\theta_1(t)-\theta_2(t)\|_{L^2(\R^2)}\|\omega_1(t)-\omega_2(t)\|_{L^2(\R^2)},
\end{split}
\end{eqnarray*}
where we have used the convention that $L^{\frac{2}{1-2\alpha}}(\R^2)$ denotes $L^\infty(\R^2)$ when $\alpha=\frac{1}{2}$.

Therefore,
\begin{eqnarray}\label{3.8}
\frac{d}{dt}\|\omega_1(t)-\omega_2(t)\|_{L^2(\R^2)}\leq CM\|\theta_1(t)-\theta_2(t)\|_{L^2(\R^2)}.
\end{eqnarray}
By direct integration in \eqref{3.8}, we obtain
\begin{eqnarray*}\begin{split}
&\|\omega_1(t)-\omega_2(t)\|_{L^2(\R^2)}\\ \leq&\|\omega_1(0)-\omega_2(0)\|_{L^2(\R^2)}+CM\int_{0}^{t}\|\theta_1(\tau)-\theta_2(\tau)\|_{L^2(\R^2)}d\tau\\
\leq& CM\int_{0}^{T}\|\theta_1(\tau)-\theta_2(\tau)\|_{L^2(\R^2)}d\tau\\
\leq &CMT\|\theta_1-\theta_2\|_{C([0,T];L^2(\R^2))}
\end{split}
\end{eqnarray*}
for $t\in[0,T].$

Since $\omega_1-\omega_2\in C([0,T];L^2(\R^2))$, we have
\begin{eqnarray}\label{3.9}
\|\omega_1-\omega_2\|_{C([0,T];L^2(\R^2))}\leq CMT\|\theta_1-\theta_2\|_{C([0,T];L^2(\R^2))}
\end{eqnarray}
From \eqref{3.9}, we can choose sufficiently small $T>0$ such that $0<CMT<1$, therefore the mapping $\mathcal{T}: B_T\rightarrow B_T$ is a contraction mapping in the topology of $C([0,T];L^2(\R^2))$.\\[2mm]
Finally, combining step 1, step 2 and step 3, by virtue of Lemma \ref{contration}, we can conclude the proof of Lemma \ref{existence of fixed-point}.

Now we are ready to prove Theorem \ref{the-1}.

\textbf{Proof of Theorem \ref{the-1}.} By virtue of Lemma \ref{existence of fixed-point}, there exists exactly one $\theta\in B_T$ such that $\omega=\mathcal{T}\theta=\theta,$ which implies that
\begin{equation*}
\left\{\ba
&\theta_{t}+u\cdot\nabla \theta = 0, ~(x,t)\in \R^{2}\times\R_{+},\\&u =\nabla^\bot(-\Delta)^{-1+\alpha}\theta,\\
&\theta(x,0)=\theta_{0}(x).\ea\ \right.
\end{equation*}
Moreover, in view of Corollary \ref{existence of scalar transport}, we have $\theta$ is in $C([0,T];H^s(\R^2)).$
This is the existence part of Theorem \ref{the-1}.
Next we turn to the uniqueness part. Suppose that $\theta_i\in C([0,T];H^s(\R^2))$, $i=1,2$ are two solutions to the equation \eqref{GSQG} with the same initial data $\theta_0\in H^s(\R^2)$. Then
\begin{equation}\label{3.10}
\left\{\ba
&\partial_t\theta_1+u_1\cdot\nabla\theta_1=0,\\
&\partial_t\theta_2+u_2\cdot\nabla\theta_2=0.\ea\ \right.
\end{equation}
Subtracting the equations in \eqref{3.10}, we obtain
\begin{eqnarray}\label{3.11}
\partial_t(\theta_1-\theta_2)+(u_1-u_2)\cdot\nabla\theta_1+u_2\cdot\nabla(\theta_1-\theta_2)=0.
\end{eqnarray}
Taking the inner product of \eqref{3.11} with $\theta_1-\theta_2$ in $L^2(\R^2)$ and using the divergence-free condition, integration by parts, H\"{o}lder inequality, Lemma \ref{CZ inequality} (used when $\alpha=\frac{1}{2}$), Corollary \ref{coro 2.3} (used when $\alpha\in(0,\frac{1}{2})$) and the standard Sobolev embedding, we can obtain
\begin{eqnarray*}
\begin{split}
&\frac{1}{2}\frac{d}{dt}\|\theta_1(t)-\theta_2(t)\|^2_{L^2(\R^2)}\\=&-\int_{\mathbb{R}^2}((u_1-u_2)\cdot\nabla\theta_1)(\theta_1-\theta_2) dx-\int_{\mathbb{R}^2}(u_2\cdot\nabla(\theta_1-\theta_2))(\theta_1-\theta_2)dx\\
=&-\int_{\mathbb{R}^2}((u_1-u_2)\cdot\nabla\theta_1)(\theta_1-\theta_2)dx\\
\leq&\|u_1-u_2\|_{L^{\frac{1}{\alpha}}(\R^2)}\|\nabla\theta_1\|_{L^{\frac{2}{1-2\alpha}}(\R^2)}\|\theta_1-\theta_2\|_{L^2(\R^2)}\\
\leq& C\|\theta_1(t)\|_{H^s(\R^2)}\|\theta_1(t)-\theta_2(t)\|^2_{L^2(\R^2)},
\end{split}
\end{eqnarray*}
where we have used the convention that $L^{\frac{2}{1-2\alpha}}(\R^2)$ denotes $L^\infty(\R^2)$ when $\alpha=\frac{1}{2}$.

Therefore,
\begin{eqnarray*}
\frac{d}{dt}\|\theta_1(t)-\theta_2(t)\|_{L^2(\R^2)}\leq C\|\theta_1(t)\|_{H^s(\R^2)}\|\theta_1(t)-\theta_2(t)\|_{L^2(\R^2)}.
\end{eqnarray*}
Applying the Gronwall lemma leads us to the following inequality,
$$\|\theta_1(t)-\theta_2(t)\|_{L^2(\R^2)}\leq e^{C\int_{0}^{t}\|\theta_1(\tau)\|_{H^s(\R^2)}d\tau}\|\theta_1(0)-\theta_2(0)\|_{L^2(\R^2)},$$ which clearly implies uniqueness of the solution.
The proof of Theorem \ref{the-1} is complete.
\section{Proof of Theorem \ref{the-2} }
In this section,  we prove Theorem \ref{the-2}.
The key difference between  $\eqref{GSQG}$ and $\eqref{tr-s-port}$ is that  the velocity in $\eqref{GSQG}$ is naturally divergence-free, while the velocity in $\eqref{tr-s-port}$ is not. However, the local well-posedness still holds true for the equation $\eqref{tr-s-port}$.

\textbf{Proof of Theorem \ref{the-2}.}
By checking the proof of  Theorem \ref{the-1}, we note that there are only two places needed to be modified to prove Theorem \ref{the-2}, that is, {\it step 3} in the proof of Lemma \ref{existence of fixed-point}  and the uniqueness part of the solution.

 To modify the {\it Step 3} in the proof of Lemma \ref{existence of fixed-point}, we can prove the contraction of the mapping $\mathcal{T}: B_T\rightarrow B_T$ as follows:
\begin{eqnarray*}
\begin{split}
&\frac{1}{2}\frac{d}{dt}\|\omega_1(t)-\omega_2(t)\|^2_{L^2(\R^2)}\\
=&-\int_{\mathbb{R}^2}((u_1-u_2)\cdot\nabla\omega_1)(\omega_1-\omega_2) dx-\int_{\mathbb{R}^2}(u_2\cdot\nabla(\omega_1-\omega_2))(\omega_1-\omega_2)dx\\
=&-\int_{\mathbb{R}^2}((u_1-u_2)\cdot\nabla\omega_1)(\omega_1-\omega_2)
dx+\frac{1}{2}\int_{\mathbb{R}^2}({\rm div} u_2)(\omega_1-\omega_2)^2dx\\
\leq&\|u_1-u_2\|_{L^{\frac{1}{\alpha}}(\R^2)}\|\nabla\omega_1\|_{L^{\frac{2}{1-2\alpha}}(\R^2)}\|\omega_1-\omega_2\|_{L^2(\R^2)}+\|\nabla u_2\|_{L^\infty(\R^2)}\|\omega_1-\omega_2\|^2_{L^2(\R^2)}\\
\leq&C\|\theta_1-\theta_2\|_{L^2(\R^2)}\|\omega_1\|_{H^s(\R^2)}\|\omega_1-\omega_2\|_{L^2(\R^2)}
+C\|\theta_2\|_{H^s(\R^2)}\|\omega_1-\omega_2\|^2_{L^2(\R^2)}\\
\leq&CM\|\theta_1(t)-\theta_2(t)\|_{L^2(\R^2)}\|\omega_1(t)-\omega_2(t)\|_{L^2(\R^2)}
+CM\|\omega_1(t)-\omega_2(t)\|^2_{L^2(\R^2)},
\end{split}
\end{eqnarray*}
where we have used the following estimate
\begin{eqnarray*}
\begin{split}
\|\nabla u_2\|_{L^\infty(\R^2)}&\leq C\|\nabla u_2\|_{H^{s-2\alpha}(\R^2)}\\ &=
C\|\nabla\nabla\Lambda^{-2+2\alpha}\theta_2\|_{H^{s-2\alpha}(\R^2)}\\&\leq
C\|\theta_2\|_{H^s(\R^2)}.
\end{split}
\end{eqnarray*}
Therefore, we obtain
\begin{eqnarray*}\begin{split}
&\frac{d}{dt}\|\omega_1(t)-\omega_2(t)\|_{L^2(\R^2)}\\ \leq &CM\|\omega_1(t)-\omega_2(t)\|_{L^2(\R^2)}+CM\|\theta_1(t)-\theta_2(t)\|_{L^2(\R^2)}.\end{split}
\end{eqnarray*}
By virtue of the Gronwall's lemma, we obtain, for $t\in[0,T]$
\begin{eqnarray*}\begin{split}
&\|\omega_1(t)-\omega_2(t)\|_{L^2(\R^2)}\\ \leq& e^{CMt}[\|\omega_1(0)-\omega_2(0)\|_{L^2(\R^2)}+CM\int_{0}^{t}\|\theta_1(\tau)-\theta_2(\tau)\|_{L^2(\R^2)}d\tau]\\
\leq&CMe^{CMT}\int_{0}^{T}\|\theta_1(\tau)-\theta_2(\tau)\|_{L^2(\R^2)}d\tau\\
\leq& CMTe^{CMT}\|\theta_1-\theta_2\|_{C([0,T];L^2(\R^2))}.\end{split}
\end{eqnarray*}

Since
$$\omega_1-\omega_2\in C([0,T];L^2(\R^2)),$$ we have
\begin{eqnarray}\label{4.1}
\|\omega_1-\omega_2\|_{C([0,T];L^2(\R^2))}\leq CMTe^{CMT}\|\theta_1-\theta_2\|_{C([0,T];L^2(\R^2))}
\end{eqnarray}
 Thanks to \eqref{4.1}, we can choose sufficiently small $T>0$ such that $0<CMTe^{CMT}<1$, therefore the mapping $\mathcal{T}: B_T\rightarrow B_T$ is a contraction mapping in the topology of $C([0,T];L^2(\R^2))$.\\[2mm]
Concerning the uniqueness of the solutions, we have
\begin{eqnarray*}
\begin{split}
&\frac{1}{2}\frac{d}{dt}\|\theta_1(t)-\theta_2(t)\|^2_{L^2(\R^2)}\\
=&-\int_{\mathbb{R}^2}((u_1-u_2)\cdot\nabla\theta_1)(\theta_1-\theta_2) dx-\int_{\mathbb{R}^2}(u_2\cdot\nabla(\theta_1-\theta_2))(\theta_1-\theta_2)dx\\
=&-\int_{\mathbb{R}^2}((u_1-u_2)\cdot\nabla\theta_1)(\theta_1-\theta_2)dx+\frac{1}{2}\int_{\mathbb{R}^2}({\rm div} u_2)(\theta_1-\theta_2)^2dx\\
\leq&\|u_1-u_2\|_{L^{\frac{1}{\alpha}}(\R^2)}\|\nabla\theta_1\|_{L^{\frac{2}{1-2\alpha}}(\R^2)}\|\theta_1-\theta_2\|_{L^2(\R^2)}+\|\nabla u_2\|_{L^\infty(\R^2)}\|\theta_1-\theta_2\|^2_{L^2(\R^2)}\\
\leq&C(\|\theta_1(t)\|_{H^s(\R^2)}+\|\theta_2(t)\|_{H^s(\R^2)})\|\theta_1(t)-\theta_2(t)\|^2_{L^2(\R^2)}.
\end{split}
\end{eqnarray*}
Therefore,
\begin{eqnarray*}
\frac{d}{dt}\|\theta_1(t)-\theta_2(t)\|_{L^2(\R^2)}\leq C(\|\theta_1(t)\|_{H^s(\R^2)}+\|\theta_2(t)\|_{H^s(\R^2)})\|\theta_1(t)-\theta_2(t)\|_{L^2(\R^2)}.
\end{eqnarray*}
By virtue of Gronwall's inequality, we  obtain
\begin{eqnarray*}
\|\theta_1(t)-\theta_2(t)\|_{L^2(\R^2)}\leq e^{C\int_{0}^{t}(\|\theta_1(\tau)\|_{H^s(\R^2)}+\|\theta_2(\tau)\|_{H^s(\R^2)})d\tau}\|\theta_1(0)-\theta_2(0)\|_{L^2(\R^2)},
\end{eqnarray*}
which  implies uniqueness of the solution.
\appendix
\section{ Littlewood-Paley theory and Besov spaces}
In this appendix, we recall some basic facts about Littlewood-Paley decomposition and  inhomogeneous Besov spaces.
For more details, it is referred to \cite{[BCD]}, \cite{[MWZ]}, \cite{[Wu]} and references
therein.

Let $(\chi,\varphi)$ be a couple of smooth functions with  values in
$[0,1]$ such that $\chi$ is supported in the ball
$\big\{\xi\in\mathbb{R}^{n}\big||\xi|\leq\frac{4}{3}\big\}$,
$\varphi$ is supported in the shell
$\big\{\xi\in\mathbb{R}^{n}\big|\frac{3}{4}\leq|\xi|\leq\frac{8}{3}\big\}$
and
\begin{align*}
    \chi(\xi)+\sum_{j\in \mathbb{N}}\varphi(2^{-j}\xi)=1\quad {\rm for \ each\ }\xi\in \mathbb{R}^{n}.
\end{align*}
For every $u\in \mathcal{S}'(\mathbb{R}^{n})$, we define the dyadic blocks as
\begin{equation*}
 \Delta_{-1}u=\chi(D)u\quad\text{and}\quad   {\Delta}_{j}u:=\varphi(2^{-j}D)u\quad {\rm for\ each\ }j\in\mathbb{N}.
\end{equation*}
We shall also use the following low-frequency cut-off:
\begin{equation*}
    {S}_{j}u:=\chi(2^{-j}D)u.
\end{equation*}
It may be  easily checked that
\begin{equation*}
    u=\sum_{j\geq-1}{\Delta}_{j}u
\end{equation*}
holds in $\mathcal{S}'(\mathbb{R}^{n})$.
\begin{definition}\label{Besov space}
For $s\in \mathbb{R}$, $(p,q)\in [1,+\infty]^{2}$ and $u\in \mathcal{S}'(\mathbb{R}^{n})$, we set
\begin{equation*}
    \norm{u}_{{B}^{s}_{p,q}(\mathbb{R}^{n})}:=
    \Big(\sum_{j\geq-1}2^{jsq}
    \norm{{\Delta}_{j}u}_{L^{p}(\mathbb{R}^{n})}^{q}\Big)^{\frac{1}{q}}
    \quad\text{if}\quad q<+\infty
\end{equation*}
and
\begin{equation*}
    \norm{u}_{{B}^{s}_{p,\infty}(\mathbb{R}^{n})}:=\sup_{j\geq-1}2^{js}\norm{{\Delta}_{j}u}_{L^{p}(\mathbb{R}^{n})}.
\end{equation*}
Then we define inhomogeneous Besov spaces as
\begin{equation*}
   {B}^{s}_{p,q}(\mathbb{R}^{n}):=\big\{u\in\mathcal{S}'(\R^{n}):\norm{u}_{{B}^{s}_{p,q}(\mathbb{R}^{n})}<+\infty\big\}.
\end{equation*}
\end{definition}

\begin{definition}\label{defapp}
For $s\in\R$, $1\leq p, q,\sigma\leq\infty$, $I=[0,T]$, the inhomogeneous space-time Besov spaces are defined as $$\widetilde{L}^{q}(I;B_{p,\sigma}^{s}(\R^{n}))=\{{u\in \mathbb{D}'(I,\mathcal{S}'(\R^{n})):\norm{u}_{\widetilde{L}^{q}(I;B_{p,\sigma}^{s}(\R^{n}))}}=
\norm{2^{js}\norm{\Delta_{j}u}_{L^q(I;L^{p}(\R^{n}))}}_{l^\sigma}<\infty\}.$$
\end{definition}
It should be remarked that the usual Sobolev spaces $H^s(\R^n)$  coincide with  Besov spaces  $B_{2,2}^s(\R^n)$.
Next, we give the embeddings of inhomogenesous Besov spaces.
\begin{lemma}[\textbf{Embeddings of Besov spaces}]\label{embeddings}
Let $s, \widetilde{s}\in\R$, and $1\leq p, q, \widetilde{p}, \widetilde{q}\leq\infty$ The following continuous embeddings hold true:

\rm(1) $B^s_{p,q}(\R^n)\hookrightarrow B^{\widetilde{s}}_{p,\widetilde{q}}(\R^n)$ whenever $\widetilde{s}<s$ or $\widetilde{s}=s$ and $\widetilde{q}\geq q$.

(2) $B^s_{p,q}(\R^n)\hookrightarrow B^{s-n(\frac{1}{p}-\frac{1}{\widetilde{p}})}_{\widetilde{p},q}(\R^n)$ whenever $\widetilde{p}\geq p$.
\end{lemma}

The following lemma is the well-known Bernstein inequality.
\begin{lemma}[\textbf{Bernstein's inequality}]\label{bern}
Let $\mathcal {B}$ be a ball of  $\R^{n}$, and $\mathcal {C}$ be a
ring of $\R^{n}$. There exists a positive constant C such that for
all integer $k\geq0$, all $1\leq a\leq b \leq\infty$ and
$u\in{L^{a}(\R^{n})}$, the following estimates are satisfied:
$$\sup_{|\alpha|=k} \|\partial^{\alpha}u\|_{L^{b}(\R^{n})}\leq
C^{k+1}\lambda^{k+n(\frac{1}{a}-\frac{1}{b})}\|u\|_{L^{a}(\R^{n})},
~~supp~\widehat{u}\subset \lambda \mathcal {B},$$
$$C^{-(k+1)}\lambda^{k}\|u\|_{L^{a}(\R^{n})}\leq \sup_{|\alpha|=k}\|\partial^{\alpha}u\|_{L^{a}(\R^{n})}\leq C^{k+1}\lambda^{k}\|u\|_{L^{a}(\R^{n})},
~~supp~\widehat{u}\subset \lambda \mathcal {C}.$$
\end{lemma}
Lastly, we review the classical Calderon-Zygmund inequality in harmonic analysis (see \cite{[S]}).
\begin{lemma}[\textbf{Calderon-Zygmund inequality}]\label{CZ inequality}
For any $f\in L^p(\R^n)$, $1<p<\infty$, we have the following inequality
$$\|\mathcal{R}_jf\|_{L^p(\R^n)}\leq C_p\|f\|_{L^p(\R^n)},$$
where $C_p$ is a constant depending only on $p$, which is independent of the function $f$.
Here $\mathcal{R}_j$ denotes the j-th componet of the vectorial Riesz transform $$\mathcal{R}=(\mathcal{R}_1,\mathcal{R}_2,...,\mathcal{R}_n)=(\frac{\partial_{x_1}}{\sqrt{-\Delta}},\frac{\partial_{x_2}}{\sqrt{-\Delta}}
,...,\frac{\partial_{x_n}}{\sqrt{-\Delta}}).$$
\end{lemma}
We also remark here that the Riesz operator $\mathcal{R}$ is bounded on Sobolev spaces $H^s(\R^n)$, since its Fourier symbol $\frac{i\xi}{|\xi|}$ is bounded.

{\bf Acknowledgements.}
 The authors would like to thank Professor Quansen Jiu for his valuable discussions and constant encouragements. H. Yu was partially supported
by the National Natural Science Foundation of China (NNSFC) (No. 11901040), Beijing Natural Science Foundation (BNSF) (No. 1204030)
and   Beijing Municipal Education Commission (KM202011232020).


\end{document}